\documentclass[10pt]{article}
\usepackage{amssymb,latexsym, amsmath, enumerate}
\usepackage{color}
\usepackage{hyperref}
\hypersetup{
                CJKbookmarks=true
	colorlinks=true,
	linkcolor=blue,
	filecolor=blue,      
	urlcolor=blue,
	citecolor=blue,
}

\def\X{{\mathrm{X}}}
\def\Y{{\mathrm{Y}}}
\def\R{\mathbb R}
\def\N{\mathbb N}
\def\C{\mathbb C}
\def\H{\mathbb H}

\newtheorem{theorem}{Theorem}
\newtheorem{lemma}{Lemma}
\newtheorem{remark}{Remark}

\title{A note on gradient estimates for the heat semigroup on nonisotropic Heisenberg groups}
\author{Ye Zhang}
\date{}
\begin{document}

\renewcommand{\theequation}{\thesection.\arabic{equation}}
\setcounter{equation}{0} \maketitle

\vspace{-1.0cm}

\bigskip

{\bf Abstract.}  In this note we obtain gradient estimates for the heat semigroup on nonisotropic Heisenberg groups.
More precisely, our aim is to get the H.-Q. Li inequality on nonisotropic Heisenberg groups, which is a
generalization of the original result on the classical Heisenberg group of dimension 3 and a counterpart of
the inequality on the H-type groups. Our proof is  based on a Cheeger type inequality, which is an approach proposed by Bakry et al.

\medskip

{\bf Mathematics Subject Classification (2000):} {\bf  35B45, 35K05, 35K08, 43A85, 58J35}

\medskip

{\bf Key words and phrases:}  Gradient estimate; Heat kernel; Heat semigroup; Nonisotropic Heisenberg group

\medskip

\renewcommand{\theequation}{\thesection.\arabic{equation}}
\section{Introduction}
\setcounter{equation}{0}
The Gaussian bounds and gradient estimates for the heat kernel on nilpotent Lie groups are important in the study of properties of the associated subelliptic operators. For example, we can obtain
\begin{align}\label{grad}
|\nabla e^{h \Delta} f(g)| \le C(\Delta,G) e^{h \Delta}(|\nabla f|)(g), \quad \forall f \in C_o^\infty(G), h > 0, g \in G
\end{align}
where $e^{h \Delta}$ is the heat semigroup corresponding to the sub-Laplacian $\Delta$ and $\nabla$ is the corresponding left-invariant sub-gradient on the nilpotent Lie group $G$ if we are given estimates of the following type:
\begin{align}\label{est}
C_1(\Delta,G) Q_h(g) e^{- \frac{d^2(g)}{4h}} \le p_h(g) \le C_2(\Delta,G) Q_h(g) e^{- \frac{d^2(g)}{4h}}, \qquad |\nabla \ln(p_h(g))| \le C_3(\Delta,G) \frac{d(g)}{h},
\end{align}
where $p_h(g)$ denotes the heat kernel associated with the sub-Laplacian $\Delta$ and $d(g) = d(g,o)$ is the Carnot-Carath\'{e}odory distance between $g$ and $o$ on this group. For the inequality \eqref{grad} on the classical Heisenberg group of dimension 3, one can refer to \cite{BBBC08,Li06}. A nature generalization of the classical Heisenberg group of dimension 3 is the higher-dimensional Heisenberg groups, or isotropic Heisenberg groups. However, there are more generalized versions of the classical Heisenberg group of dimension 3, that is, nonisotropic Heisenberg groups and H-type groups. The corresponding result on H-type groups can be found in \cite{NE10,HLi10}.

In this note we give the proof of \eqref{grad} on  a nonisotropic Heisenberg group $\H = \H(\mathcal{K},\mathcal{A})$ (see Section 2) via a Cheeger type inequality, which is an approach proposed by Bakry et al. However, this approach needs an additional estimate of the heat kernel $p_h$ associated with the canonical sub-Laplacian $\Delta$ on $\H$: for $g = (z,t) \in \H$,
\begin{align}
|\partial_{t} \log{p_h(z,t)}| \le \frac{C_4(\Delta, \H)}{h}.
\end{align}
We obtain the following main theorem.

\vskip4pt

\begin{theorem}\label{t1}
There exists a constant $K = K(\H)$ such that for all $f \in C_o^\infty(\H)$ we have
\begin{align}
|\nabla e^{h \Delta} f(g)| \le K e^{h \Delta}(|\nabla f|)(g), \quad \forall h > 0, g \in \H.
\end{align}
\end{theorem}

\vskip4pt

From Theorem \ref{t1} and H\"{o}lder's inequality, we have
\begin{align}
|\nabla e^{h \Delta} f(g)| \le K (e^{h \Delta}(|\nabla f|^p))^{\frac{1}{p}}(g), \quad \forall h > 0, g \in \H
\end{align}
immediately, which is a generalization the result of \cite{DM05}.

Furthermore, from the argument of \cite[Subsection 6.1]{BBBC08}, we have the following
logarithmic Sobolev inequality and Poincar\'{e} inequality:
\begin{align}\label{lse}
e^{h \Delta}( f^2\log{f^2}) - e^{h \Delta}(f^2) \log(e^{h \Delta}(f^2)) &\le K^{\prime} h e^{h \Delta}(|\nabla f|^2),\\
\label{poe}
e^{h \Delta}( f^2) - (e^{h \Delta}f)^2 &\le K^{\prime\prime} h e^{h \Delta}(|\nabla f|^2), \quad \forall \, f \in C_o^{\infty}(\H),
\end{align}
where $K^{\prime} = K^{\prime}(\H)$ and $K^{\prime\prime} = K^{\prime\prime}(\H)$ are two constants depending on the group data $\H$. In fact, as shown in \cite{GL21}, the constants $K^{\prime}$ and $K^{\prime\prime}$ can be chosen as universal constants (which are independent of $\H$). Their proof bases on the inequality \eqref{lse} on the classical Heisenberg group of dimension 3 and the properties of the logarithmic Sobolev inequality under tensorization and quotient, which is different from our approach.

The note is organized as follows. In Section \ref{st2} we give a review of nonisotropic Heisenberg groups.
In Section \ref{st3} we establish a Cheeger type inequality on the nonisotropic Heisenberg group. This can be done by introducing ``polar coordinates''
which appear in \cite{BBBC08} and simplified by Eldredge in \cite{NE10}.  The proof of Theorem \ref{t1} will be given
in Section \ref{st4}.

Our proof is adapted from the proof in \cite{BBBC08} and \cite{NE10}, so we only sketch the proof when there are few changes.

\section{Preliminaries} \label{st2}
\setcounter{equation}{0}
\subsection{The nonisotropic Heisenberg group}
Suppose that
\begin{align*}
l \in \N^+, \  \mathcal{K} = (k_1, \cdots, k_l) \in \left( \N^+ \right)^l,  \  \mathcal{A} = (a_1, \cdots, a_l)  \mbox{ with } 0 <  a_1 < \cdots < a_l = 1, \ \mbox{and} \  n = \sum_{i = 1}^l k_i.
\end{align*}
Recall that (see \cite{BGG00} or \cite{LP03}) the general Heisenberg group with parameter $(\mathcal{K}, \mathcal{A})$,  $\H(\mathcal{K}, \mathcal{A}) = \left( \prod\limits_{i = 1}^l \C^{k_i} \right) \times \R$, is a stratified group
with the group multiplication
\begin{align*}
(z, t) \cdot (z', t') = \left(z + z', t + t' + 2 \sum_{i = 1}^l a_i \Im
\langle z_i, z_i' \rangle\right),
\end{align*}
where $z = (z_1, \cdots, z_l)$, $z' =
(z_1', \cdots, z_l') \in \prod\limits_{i = 1}^l \C^{k_i}$, $\Im \omega$ denotes the imaginary part of a complex number $\omega$ and $\langle \cdot, \cdot \rangle$ denotes the usual complex inner product. When $l = 1$, $\H(n, 1)$ is often called the Heisenberg group or the isotropic Heisenberg group of real dimension $2 n + 1$. When $l \geq 2$, $\H(\mathcal{K}, \mathcal{A})$ is often called the non-isotropic Heisenberg group of real dimension $2 n + 1$. For convenience, we sometimes drop the depending parameter $(\mathcal{K}, \mathcal{A})$ and denote the group by $\H$.
Recall that the dilation $\delta_r$ ($r > 0$) is given by
\begin{align*}
\delta_{r}(z, t) = (r z, r^2 t), \quad \forall (z, t) \in \H(\mathcal{K}, \mathcal{A}) = \C^n \times \R.
\end{align*}

For $1 \leq i \leq l$, set
\begin{align*}
z_i = (z_{i, 1}, \cdots, z_{i, k_i}) = (x_{i, 1} + \imath y_{i, 1}, \cdots, x_{i, k_i} + \imath y_{i, k_i})
\end{align*}
with $x_{i, j}, y_{i, j} \in \R, 1 \leq j \leq k_i$,
$|z_i| = \sqrt{\langle z_i, z_i \rangle}$ and $|z| = \sqrt{\sum\limits_{i = 1}^l |z_i|^2}$.
Let
\begin{align*}
\X_{i, j} &= \frac{\partial}{\partial x_{i, j}} + 2 a_i y_{i, j} \frac{\partial}{\partial t},  \quad \Y_{i, j} = \frac{\partial}{\partial
y_{i, j}} - 2 a_i x_{i, j} \frac{\partial}{\partial t}, \\
\hat{\X}_{i, j} &= \frac{\partial}{\partial x_{i, j}} - 2 a_i y_{i, j} \frac{\partial}{\partial t},  \quad \hat{\Y}_{i, j} = \frac{\partial}{\partial
y_{i, j}} + 2 a_i x_{i, j} \frac{\partial}{\partial t}, \qquad 1 \leq i \leq l, 1 \leq j \leq k_i,
\end{align*}
be the left and right invariant vector fields respectively on $\H(\mathcal{K}, \mathcal{A})$. The associated sub-Laplacian, left-invariant sub-gradient and
right-invariant sub-gradient are given respectively by
\begin{align*}
\Delta = \sum_{i = 1}^l \sum_{j = 1}^{k_i} (\X_{i, j}^2 + \Y_{i, j}^2), \qquad \nabla = (\X_{1, 1}, \Y_{1, 1}, \cdots, \X_{l, k_l}, \Y_{l, k_l}),
\qquad \hat{\nabla} = (\hat{\X}_{1, 1}, \hat{\Y}_{1, 1}, \cdots, \hat{\X}_{l, k_l}, \hat{\Y}_{l, k_l}).
\end{align*}
Let $d(g,g^{\prime})$ denote the Carnot-Carath\'eodory distance from $g$ to $g^{\prime}$, and $p_h$ ($h
> 0$) the heat kernel (that is, the integral kernel of $e^{h \Delta}$)
on $\H(\mathcal{K}, \mathcal{A})$. Let $o = (0, \cdots, 0)$ be the origin of $\H(\mathcal{K}, \mathcal{A})$,
$g = (z, t) \in \C^n \times \R$. We set in the sequel 
\begin{align*}
d = d(g) = d(g, o), \qquad p_h(g) = p_h(g, o) \quad \mbox{and} \quad p(g) = p_1(g).
\end{align*}

It is well known that 
\begin{align} \label{sp}
p_h(z, t) = \frac{1}{h^{n+1}} p\left(\frac{z}{\sqrt{h}}, \frac{t}{h}\right), \qquad \forall h > 0, \quad (z, t) \in \C^n \times \R,
\end{align}
and   
(cf. \cite{BGG00} or \cite{LP03}):
\begin{equation} \label{1}
p_h(z, t)  
= \frac{1}{2 (4\pi h)^{n+1}} \int_{\R}\prod_{j=1}^l \left( \frac{a_j \lambda}{\sinh{(a_j \lambda)}} \right)^{k_j} \exp{\left( \frac{1}{4 h} (\imath \lambda t - \sum_{j=1}^l |z_j|^2 a_j \lambda \coth{(a_j \lambda)}) \right)} \, d\lambda.
\end{equation}

We notice that the Lebesgue measure $m$ on $\C^n \times \R $  is invariant under the group operation (including left, right multiplication and inversion) and consequently it is the Haar measure on
$\H(\mathcal{K},\mathcal{A})$. Set $e^{h \Delta}$ be the heat semigroup corresponding to the sub-Laplacian $\Delta$, we have (see \cite[Lemma 2]{LP03})
\begin{align}
e^{h\Delta}f(g) = f \ast p_h (g)  & = \int_{\H} f(g^{\prime}) p_h((g^{\prime})^{-1} \cdot g) m(dg^{\prime})\\
\label{semigroup}
&=\int_{\H} f(g \cdot g^{\prime}) p_h(g^{\prime}) m(dg^{\prime}), \qquad \forall h>0, f \in C_o^{\infty}(\H)
\quad \mbox{and} \quad g \in \H(\mathcal{K}, \mathcal{A}),
\end{align}
where we use the fact that $p_h(g) = p_h(g^{-1})$ from the expression of the heat kernel (cf. \eqref{1}) and the invariance of the Lebesgue measure under the group operation.

From \eqref{semigroup} we can obtain that for $1 \le i \le l, 1 \le j \le k_i$,
\begin{align}\label{rightcom}
\hat{\X}_{i,j}e^{h\Delta} = e^{h\Delta}\hat{\X}_{i,j}, \quad \hat{\Y}_{i,j}e^{h\Delta} = e^{h\Delta}\hat{\Y}_{i,j}.
\end{align}

\subsection{The Carnot-Carath\'{e}odory distance on $\H(\mathcal{K}, \mathcal{A})$}

We first recall the Carnot-Carath\'{e}odory distance on $\H(\mathcal{K}, \mathcal{A})$ (cf. \cite[\S 3]{BGG00}). Set
\begin{align} \label{mu}
\mu(\omega) = \frac{\omega}{{\sin^2{\omega}}}-\cot{\omega} = \frac{2 \omega - \sin{2 \omega}}{2 \sin^2{\omega}}: \quad (-\pi, \pi) \longrightarrow \R,
\end{align}
which is a monotonely increasing diffeomorphism between $(-\pi,\pi)$ and $\R$. Given $(z,t)$ such that $z_l \ne 0$, there exists exactly one $-\pi < \theta = \theta(z, t) < \pi$ such that
\begin{align} \label{dd1}
t = \sum_{j=1}^l a_j \mu(a_j \theta) |z_j|^2,
\end{align}
and we have
\begin{align} \label{dd2}
d^2(z,t) = \sum_{j=1}^l \left( \frac{a_j \theta}{\sin(a_j \theta)} \right)^2 |z_j|^2
= \theta \, ( t + \sum_{j=1}^l a_j \cot{(a_j\theta)} \, |z_j|^2 ).
\end{align}

In the other case $z_l = 0$, we have the following alternatives:\\
1. If $|t| < \sum\limits_{j=1}^{l - 1} a_j \mu(a_j \pi) |z_j|^2$,  the equation (\ref{dd1}) has a unique solution $-\pi < \theta < \pi$ and \eqref{dd2} is still valid. \\ 
2. If $|t| \geq \sum\limits_{j=1}^{l - 1} a_j \mu(a_j \pi) |z_j|^2$, then we have
\begin{align*}
d^2(z, t) =
\pi \left( |t| + \sum_{j=1}^{l-1} a_j \cot{(a_j \pi)} |z_j|^2 \right).
\end{align*}

Using the equivalence between the Carnot-Carath\'eodory distance and a homogeneous norm on stratified groups (see for example \cite{VSC92}), or by a
direct calculation, we have (see the Subsection \ref{st23} below for the notation ``$\sim$'')
\begin{align} \label{l1}
d^2(z, t) \sim |z|^2 + |t|, \qquad \forall (z, t) \in \H(\mathcal{K}, \mathcal{A}).
\end{align}

\medskip

\subsection{Notation}\label{st23}
Let $f$ and $w$ be two real-valued functions. From now on, the notation $f \lesssim w$ means that there exists a universal constant $C > 0$, maybe depending on $(\mathcal{K}, \mathcal{A})$, such that $f \leq C w$. We use the notation $f \sim w$ if $f \lesssim w$ and $w \lesssim f$. Furthermore, $f \lesssim_{\vartheta} w$, with a parameter $\vartheta$, means that there exists a constant $C(\vartheta) > 0$, depending on $\vartheta$ and maybe depending on $(\mathcal{K}, \mathcal{A})$, such that $f \leq C(\vartheta) w$. Similarly, we use the notation $f \sim_{\vartheta} w$ if $f \lesssim_{\vartheta} w$ and $w \lesssim_{\vartheta} f$.

\vskip4pt

Let $A_1$ and $A_2$ be two measurable sets of $\H(\mathcal{K}, \mathcal{A})$. We use the notation $A_1 \doteq A_2$ if their
symmetric difference $A_1 \triangle A_2$ has Lebesgue measure zero.

\subsection{Some results on $\H(\mathcal{K}, \mathcal{A})$}
The results in this subsection play a crucial role in our proof of Theorem \ref{t1}.

\vskip4pt

These following two lemmas can be found in \cite[Corollary 2 and Theorem 5]{LZ18}.

\vskip4pt

\begin{lemma}\label{l1}
We have
\begin{align} \label{pehk}
p(z, t) \sim \left( 1 + |z|^2 \epsilon_0^2 + |z_l|^2 \epsilon_0^{-1} \right)^{- \frac{1}{2}}
 \left( \frac{1 + |z| + |z_l|^2  \epsilon_0^{-2}}{1 + |z| \epsilon_0 + |z_l|^2 \epsilon_0^{-1}} \right)^{k_l - 1} e^{-\frac{d^2(z, t)}{4}},
\end{align}
for all $(z, t) = ((z_1, \cdots, z_l), t) \in \H(\mathcal{K}, \mathcal{A})$ satisfying $\epsilon_0 = \frac{\sin{\theta}}{\theta} > 0$. \\
\end{lemma}

\vskip4pt

\begin{lemma}\label{l2} We have
$| \nabla \ln{p(z,t)} | \lesssim d(z,t)$ and $| \partial_t\ln{p(z,t)} | \lesssim 1 $ for all $(z,t) \in \H(\mathcal{K}, \mathcal{A})$.
\end{lemma}

\vskip4pt

\begin{remark}
In fact, \cite[Theorem 5]{LZ18} only states the first part of Lemma \ref{l2}. However, the second part of Lemma \ref{l2} follows easily from the proof of \cite[Theorem 5]{LZ18}.
\end{remark}

\vskip4pt

Let $B$ be the  Carnot-Carath\'{e}odory unit ball, that is, $B=\{g|d(g) < 1\}$ and for $f \in C_o^\infty(\H)$ we set $m_f = \frac{\int_B f dm}{\int_B dm}$. Then we have
the following Poincar\'{e} inequality, which can be found in \cite{MS95}.

\vskip4pt

\begin{lemma}\label{l3}
 For any $f \in C_o^{\infty}(\H)$, we have
\begin{align}\label{Poincare}
\int_B |f - m_f| dm \lesssim  \int_B |\nabla f| dm.
\end{align}
\end{lemma}

\section{A Cheeger type inequality} \label{st3}
\setcounter{equation}{0}
In this section we establish the following theorem, which is a Cheeger type inequality (see \cite[Theorem 6.3]{BBBC08}).

\vskip4pt

\begin{theorem}\label{t2}
For all $f \in C_o^\infty(\H)$, we have
\begin{align}
\int_{\H} |f - m_f| p dm \lesssim \int_{\H} |\nabla f| p dm.
\end{align}
\end{theorem}

\vskip4pt

From Lemma \ref{l3} and the fact that $p$ is positive and has a upper and lower bound on $B$, the main difficulty occurs in the integral on $B^c$ and the following lemma deals with it.
As a result, Theorem \ref{t2} follows from these two lemmas.

\vskip4pt

\begin{lemma}\label{l4}
For all $f \in C_o^\infty(\H)$, we have
\begin{align}\label{Poincare2}
\int_{B^c} |f - m_f| p dm \lesssim \int_{\H} |\nabla f| p dm.
\end{align}
\end{lemma}

\vskip4pt

To prove this lemma, we introduce a change of variables (so-called ``polar coordinates'') $\Psi: \{(u,\eta)\in \prod\limits_{j=1}^l\C^{k_j} \times \R | u_l \ne 0, 0< |\eta| < \pi \} \longrightarrow
\{(z,t) \in \H(\mathcal{K}, \mathcal{A})| z_l \ne 0, t \ne 0 \}$, where
\begin{align}\label{psi}
\Psi(u,\eta) = \left( (1 - e^{- 2\imath a_1 \eta})u_1, \cdots,  (1 - e^{-2 \imath a_l \eta})u_l, \sum_{j=1}^l 2a_j |u_j|^2(2a_j \eta - \sin(2a_j \eta)) \right).
\end{align}
In this new coordinates, set $U = \left(4\sum\limits_{j=1}^l a_j^2 |u_j|^2 \right)^{\frac{1}{2}}$  and we have
\begin{align*}
|z_j|^2 &= |u_j|^2 (2 - 2 \cos(2a_j\eta)), \quad \theta = \eta \in (-\pi, 0) \cup (0 , \pi), \quad d^2(u,\eta) = d^2(\Psi(u,\eta)) = 4 \sum_{j=1}^l
a_j^2 |u_j|^2 \eta^2 = U^2 \eta^2 \\
B &\doteq \left\{(u,\eta)\Big{|} 0 < |\eta| <\min\left(\frac{1}{U},\pi\right)\right\}, \quad \mbox{and} \quad
B^c \doteq \left\{(u,\eta)\Big{|} \frac{1}{U} \le |\eta| < \pi \right\}.
\end{align*}
Here ``$\doteq$'' means that the right hand side is in fact smaller than the left hand side and they differ on a set of measure zero
(see Subsection \ref{st23}).
As a consequence, it will not affect the value of the integral. \\
We set $f(u, \eta) = f(\Psi(u, \eta))$, $p(u,\eta) = p(\Psi(u,\eta))$ and $u^{\prime} = (u_1 , \cdots, u_{l-1}) \in \prod\limits_{j=1}^{l-1} \C^{k_j}$ and we have
\begin{align}
\nonumber
\frac{d}{ds} f(u,s\eta) &= 2\eta \sum_{i = 1}^l \Big[ a_i \sum_{j = 1}^{k_i} (\sin(2a_i s \eta)\Re{u_{i,j}} - \cos(2 a_i s \eta) \Im{u_{i,j}}) \X_{i,j}f(u,s\eta)\\
{\label{horizon}}
&+(\sin(2 a_i s \eta) \Im{u_{i,j}} + \cos(2 a_i s \eta)\Re{u_{i,j}}) \Y_{i,j}f(u,s\eta) \Big],
\end{align}
where $\Re \omega$ denotes the real part of a complex number $\omega$.
From \eqref{horizon} we know that for fixed $(u,\eta)$, $s \to \Psi(u,s\eta)$ is the shortest horizontal path from the origin to
$\Psi(u,\eta)$ (cf. \cite[Section 3]{NE09} or \cite[Subsection 2.1]{ABGR09}).
Furthermore, \eqref{horizon} and Cauchy-Schwarz inequality yield
\begin{align}\label{dds}
\left|\frac{d}{ds}f(u,s\eta)\right| \le  U|\eta||\nabla f(u,s\eta)|.
\end{align}
Under the new coordinates, \eqref{pehk} becomes
\begin{align}\label{npehk}
p(u,\eta) \sim  \left\{\begin{array}{ll}
1 & \textrm{if $U|\eta| \le 1$}\\
\frac{1}{|u||\eta|}e^{-\frac{U^2\eta^2}{4}} & \textrm{if $U|\eta| \ge 1$ and $\pi - |\eta| \ge \frac{\pi}{8}$}\\
\frac{1}{(|u^{\prime}|^2(\pi - |\eta|) +|u_l|^2 )^{\frac{1}{2}}(\pi - |\eta|)^{k_l-\frac{1}{2}}}e^{-\frac{U^2\eta^2}{4}} & \textrm{if $U|\eta| \ge 1$, $\pi - |\eta| \le \frac{\pi}{4}$   }\\
 & \textrm{and $|u^{\prime}|^2(\pi - |\eta|)^2 +|u_l|^2(\pi - |\eta|) \ge 100$} \\
\left(|u_l|^2 + |u^{\prime}| + |u_l|(\pi - |\eta|)\right)^{k_l-1}e^{-\frac{U^2\eta^2}{4}} & \textrm{if $U|\eta| \ge 1$, $\pi - |\eta| \le \frac{\pi}{4}$ }\\
 & \textrm{and $|u^{\prime}|^2(\pi - |\eta|)^2 +|u_l|^2(\pi - |\eta|) \le 100$}
\end{array} \right. .
\end{align}

We denote the Jacobian determinant of $\Psi$ by $J(u,\eta)$ and we are in a place to compute it. We begin with a lemma, which is easy to prove.

\vskip4pt

\begin{lemma}\label{l5}
For $m \ge 4$ and
\begin{align*}
M_{m \times m} =
\begin{bmatrix}
b_1 & b_2 & 0 & \cdots & 0 & b_3 \\
b_4 & b_5 & 0 & \cdots & 0 & b_6 \\
 0  &   0 &   &  &   &  \\
\vdots&\vdots &  &\text{{\huge $Q$}} & & \\
 0 & 0 & & & &  \\
b_7 & b_8 &  &  &  &
\end{bmatrix}_{m \times m}
\mbox{with} \quad Q_{(m - 2) \times (m - 2)} =
\begin{bmatrix}
 & & & \ast \\
 &\text{{\huge $Q_1$}} & & \ast \\
 & & & \ast \\
 \ast& \ast & \ast & \ast
\end{bmatrix}_{(m -2) \times (m -2)},
\end{align*}
we have
\begin{align}
|M| = (b_1 b_5 - b_2 b_4) |Q| + (b_3 b_4 b_8 + b_2 b_6 b_7 - b_1 b_6 b_8 - b_3 b_5 b_7)|Q_1|.
\end{align}
\end{lemma}

\vskip4pt

From direct computation we have the Jacobian matrix is
\begin{align}
\begin{bmatrix}
1-\cos(2a_1\eta) & -\sin(2a_1\eta) & 0 & \cdots & 0 & 0 & v_1 \\
\sin(2a_1\eta)  & 1-\cos(2a_1\eta) & 0 &  \cdots & 0 &0  & v_2 \\
  0 & 0 & & & & & \\
 \vdots & \vdots &  &\ddots& & & \vdots \\
  0 & 0 & & & 1-\cos(2a_l\eta)  &-\sin(2a_l\eta)  & v_{2n -1}\\
 0 & 0 &  & & \sin(2a_l\eta) &1-\cos(2a_l\eta) & v_{2n}\\
 r_1 & r_2 &  & \cdots & r_{2n-1} & r_{2n} &\sum\limits_{j=1}^l 4a_j^2 |u_j|^2 (1-\cos(2a_j \eta))
\end{bmatrix},
\end{align}
where
\begin{align*}
\begin{bmatrix}
v_1 \\
v_2 \\
\vdots \\
v_{2n-1} \\
v_{2n}
\end{bmatrix}
=
\begin{bmatrix}
2a_1 \sin(2a_1\eta)\Re u_{1,1} - 2a_1 \cos(2a_1\eta) \Im u_{1,1} \\
2a_1 \sin(2a_1 \eta) \Im u_{1,1} + 2a_1 \cos(2a_1 \eta) \Re u_{1,1} \\
\vdots \\
2a_l \sin(2a_l\eta)\Re u_{l,k_l} - 2a_l \cos(2a_l\eta) \Im u_{l,k_l} \\
2a_l \sin(2a_l \eta) \Im u_{l,k_l} + 2a_l \cos(2a_l \eta) \Re u_{l,k_l}
\end{bmatrix}
\mbox{and} \quad
\begin{bmatrix}
r_1 \\
r_2 \\
\vdots \\
r_{2n-1} \\
r_{2n}
\end{bmatrix}^{\mathrm{T}}
=
\begin{bmatrix}
4a_1 \Re u_{1,1}(2a_1 \eta - \sin(2a_1 \eta))\\
4a_1 \Im u_{1,1}(2a_1 \eta - \sin(2a_1 \eta))\\
\vdots \\
4a_l \Re u_{l,k_l}(2a_l \eta - \sin(2a_l \eta))\\
4a_l \Im u_{l,k_l}(2a_l \eta - \sin(2a_l \eta))
\end{bmatrix}^{\mathrm{T}}.
\end{align*}
To compute the determinant, we use induction and the Lemma \ref{l5} above. The result is
\begin{align}
J(u,\eta) = \sum_{j=1}^l \left[ \prod_{i \ne j} (2 - 2\cos(2a_i \eta))^{k_i}\right] 8 a_j^2 |u_j|^2 (2-2 \cos(2a_j \eta) -
2a_j \eta \sin(2a_j \eta)) (2 - 2\cos(2a_j \eta))^{k_j-1}.
\end{align}
As a result, we have
\begin{align}\label{Jacobi}
J(u,\eta)\sim |\eta|^{2n+2}(\pi - |\eta|)^{2k_l -1}\left(|u^{\prime}|^2(\pi- |\eta|)+ |u_l|^2\right)
\end{align}
and for $(u,\eta) \in B^c$ we have
\begin{align}
p(u,\eta)J(u,\eta) \sim \left\{\begin{array}{ll}
|u||\eta|^{2n+1}e^{-\frac{U^2\eta^2}{4}} & \textrm{if $\pi - |\eta| \ge  \frac{\pi}{8}$}\\
(|u^{\prime}|^2(\pi - |\eta|) +|u_l|^2 )^{\frac{1}{2}}(\pi - |\eta|)^{k_l-\frac{1}{2}}e^{-\frac{U^2\eta^2}{4}} &
\textrm{if $\pi - |\eta| \le \frac{\pi}{4}$ and }\\
 & \textrm{$|u^{\prime}|^2(\pi - |\eta|)^2 +|u_l|^2(\pi - |\eta|) \ge 100$} \\
\left(|u_l|^2 + |u^{\prime}| + |u_l|(\pi - |\eta|)\right)^{k_l-1} (\pi - |\eta|)^{2k_l -1}
& \textrm{if  $\pi - |\eta| \le \frac{\pi}{4}$ and}\\
\left(|u^{\prime}|^2(\pi- |\eta|)+ |u_l|^2\right) e^{-\frac{U^2\eta^2}{4}} & \textrm{ $|u^{\prime}|^2(\pi - |\eta|)^2 +|u_l|^2(\pi - |\eta|) \le 100$}
\end{array} \right. .
\end{align}
Then we establish the following lemma and with this lemma we can follow the proof of \cite[Lemma 3.9]{NE10} to obtain Lemma \ref{l4}.

\vskip4pt

\begin{lemma}\label{l6}
If $(u,\eta) \in B^c$ we have
\begin{align}
\int_{1}^{\frac{\pi}{|\eta|}} p(u,v\eta)J(u,v\eta) dv \lesssim \frac{1}{|u|^2|\eta|^2} p(u,\eta)J(v,\eta).
\end{align}
\end{lemma}

\vskip4pt

\begin{remark}
Lemma \ref{l6} is a counterpart of a special case of the result of \cite[Lemma 3.4]{NE10} (that is, $q = 0$).
However, it suffices to use Lemma \ref{l6} to prove Theorem \ref{t2} (one can check the proof of \cite[Lemma 3.9]{NE10}).
In fact, our proof here can be adapted to the general case $q \ne 0$ without difficulties.
\end{remark}

\vskip4pt

 \noindent \textbf{Proof.}
We first split the region $B^c$ into three parts. To be more precise, we pick $\theta_0 = \frac{\pi}{4}$ and $\gamma_0 = 100$. Then we set
\begin{align*}
R_1 &= \left\{(u,\eta) \in B^c | 0 < |\eta| < \theta_0 \right\} \\
R_2 &= \left\{(u,\eta) \in B^c | \theta_0 < |\eta| < \pi, |u^{\prime}|^2(\pi - |\eta|)^2 +|u_l|^2(\pi - |\eta|) > \gamma_0 \right\} \\
R_3 &= \left\{(u,\eta) \in B^c | \theta_0 < |\eta| < \pi, |u^{\prime}|^2(\pi - |\eta|)^2 +|u_l|^2(\pi - |\eta|) \le  \gamma_0 \right\}
\end{align*}
and we have $B^c = \mathop{\bigcup}\limits_{j = 1}^3 R_j$.

We first consider the case that $(u,\eta) \in R_3$. In this case, for any $v \in \left[1, \frac{\pi}{|\eta|}\right)$ we have
$(u,v\eta) \in R_3$. Noticing that from monotonicity for $v \in \left[1, \frac{\pi}{|\eta|}\right)$,
\begin{align*}
 &\quad \left(|u_l|^2 + |u^{\prime}| + |u_l|(\pi - v|\eta|)\right)^{k_l-1} (\pi - v|\eta|)^{2k_l -1}
\left(|u^{\prime}|^2(\pi- v|\eta|)+ |u_l|^2\right) \\
 &\le
\left(|u_l|^2 + |u^{\prime}| + |u_l|(\pi - |\eta|)\right)^{k_l-1} (\pi - |\eta|)^{2k_l -1}
\left(|u^{\prime}|^2(\pi- |\eta|)+ |u_l|^2\right),
\end{align*}
we have
\begin{align*}
p(u,v\eta)J(u,v\eta) \lesssim p(u,\eta)J(u,\eta) e^{-\frac{1}{4}U^2 \eta^2 (v^2 - 1)}.
\end{align*}
Consequently,
\begin{align*}
\int_{1}^{\frac{\pi}{|\eta|}} p(u,v\eta)J(u,v\eta) dv
&\lesssim p(u,\eta)J(u,\eta) \int_{1}^{\frac{\pi}{|\eta|}} e^{-\frac{1}{4}U^2 \eta^2 (v^2 - 1)} dv \\
&\lesssim p(u,\eta)J(u,\eta) \int_0^{\frac{\pi^2}{\eta^2} - 1} e^{-\frac{1}{4}U^2 \eta^2 w} \frac{1}{\sqrt{w+1}} dw \\
&\lesssim \frac{1}{|u|^2|\eta|^2} p(u,\eta)J(v,\eta)
\end{align*}

Then we consider the case that $(u,\eta) \in R_2$. In this case, there exist a $v_0 \in \left(1, \frac{\pi}{|\eta|}\right)$ such that
for $v \in (1,v_0)$ we have $(u,v\eta) \in R_2$ and for $v \in \left[v_0, \frac{\pi}{|\eta|}\right)$ we have $(u,v\eta) \in R_3$. That is, $|u^{\prime}|^2(\pi - v_0|\eta|)^2 +|u_l|^2(\pi - v_0|\eta|) = \gamma_0 = 100$. As a result, we split the integral into two parts.
\begin{align*}
\int_{1}^{\frac{\pi}{|\eta|}} p(u,v\eta)J(u,v\eta) dv  = \int_1^{v_0} + \int_{v_0}^{\frac{\pi}{|\eta|}} = I_1 + I_2.
\end{align*}
We begin by the estimate of $I_2$. As before, from monotonicity we obtain that
\begin{align*}
I_2 &\lesssim \frac{1}{|u|^2|\eta|^2} p(u,v_0\eta)J(u,v_0\eta) \\
&\lesssim \frac{1}{|u|^2|\eta|^2}(|u^{\prime}|^2(\pi - v_0|\eta|) +|u_l|^2 )^{\frac{1}{2}}(\pi - v_0|\eta|)^{k_l-\frac{1}{2}}e^{-\frac{U^2\eta^2}{4}}
\lesssim \frac{1}{|u|^2|\eta|^2} p(u,\eta)J(u,\eta).
\end{align*}
Similarly we have
\begin{align*}
I_1 \lesssim  p(u,\eta)J(u,\eta) \int_{1}^{v_0} e^{-\frac{1}{4}U^2 \eta^2 (v^2 - 1)} dv
\lesssim \frac{1}{|u|^2|\eta|^2} p(u,\eta)J(u,\eta).
\end{align*}

Finally we consider the case that $(u,\eta) \in R_1$. In this case, we pick $\theta_1 = \frac{\pi}{8}$
and split the integral into two parts.
\begin{align*}
\int_{1}^{\frac{\pi}{|\eta|}} p(u,v\eta)J(u,v\eta) dv  = \int_1^{\frac{\theta_0 + \theta_1}{|\eta|}} +
\int_{\frac{\theta_0 + \theta_1}{|\eta|}}^{\frac{\pi}{|\eta|}} = J_1 + J_2.
\end{align*}
We first estimate $J_2$. We note that
\begin{align}
p(u,\eta) \lesssim  \left(|u_l|^2 + |u^{\prime}| + |u_l|(\pi - |\eta|)\right)^{k_l-1}  e^{-\frac{U^2\eta^2}{4}}
\end{align}
when $\pi - |\eta| \le \frac{\pi}{4}$ and $|u^{\prime}|^2(\pi - |\eta|)^2 +|u_l|^2(\pi - |\eta|) \ge 100$.
As a result, we only need to prove for any $r \ge 0$ we have
\begin{align*}
\int_{\frac{\theta_0 + \theta_1}{|\eta|}}^{\frac{\pi}{|\eta|}} |u|^{r} e^{-\frac{1}{4}U^2\eta^2v^2} dv \lesssim_r |\eta|^{2n - 1}e^{-\frac{U^2\eta^2}{4}}
\end{align*}
In fact, since $v \in \left[\frac{\theta_0 + \theta_1}{|\eta|}, \frac{\pi}{|\eta|}\right)$, we have
$\eta^2(v^2 -1) \ge \eta^2 \frac{\theta_1}{|\eta|}\frac{\theta_1}{|\eta|} = \theta_1^2$. Consequently we have
$|u|^{r+2n}e^{-\frac{1}{4}U^2\eta^2(v^2-1)} \lesssim_r 1$ and
\begin{align*}
\int_{\frac{\theta_0 + \theta_1}{|\eta|}}^{\frac{\pi}{|\eta|}} |u|^{r}  e^{-\frac{1}{4}U^2\eta^2(v^2-1)} dv
\lesssim_r \frac{1}{|u|^{2n}} \frac{1}{|\eta|} \lesssim |\eta|^{2n - 1}
\end{align*}
since $|u||\eta| \sim U|\eta| \ge 1$, which gives our desired estimate. Then we are in a place to estimate $J_1$. However, \cite[Lemma 3.3]{NE10} yields
\begin{align*}
J_1 \lesssim |u||\eta|^{2n+1} \int_1^{+\infty} v^{2n+1} e^{-\frac{1}{4}U^2\eta^2v^2} dv
\lesssim \frac{|u||\eta|^{2n+1}}{|u|^2|\eta|^2} e^{-\frac{1}{4}U^2\eta^2} \sim \frac{1}{|u|^2|\eta|^2} p(u,\eta)J(u,\eta),
\end{align*}
which ends our proof in this case. Combining all results above, we proves Lemma \ref{l6}.
~ \hspace*{20pt} ~ \hfill $\Box$

\section{Proof of Theorem 1} \label{st4}
\setcounter{equation}{0}
First, from translation and dilation properties we only need to prove that there exists a constant $K$ such that for all $f \in C_o^\infty(\H)$
we have
\begin{align}\label{reduce}
|\nabla e^{ \Delta} f(0)| \le K e^{ \Delta}(|\nabla f|)(0).
\end{align}
From the Markov property $e^{h\Delta}(1) = 1$ we can further assume that $m_f = 0$.
Noticing that $\hat{\nabla}$ is commute with $e^{h \Delta}$ (cf. \eqref{rightcom}) and $\nabla = \hat{\nabla}$ at $0$, we have
\begin{align*}
\nabla e^{ \Delta} f(0) = \hat{\nabla}e^{ \Delta} f(0) = e^{ \Delta} \hat{\nabla} f(0).
\end{align*}
Then (\ref{reduce}) is reduced to
\begin{align*}
|e^{ \Delta} \hat{\nabla} f(0)| \le K e^{ \Delta}(|\nabla f|)(0)
\end{align*}
which equivalently is
\begin{align*}
\left| \int_{\H} \hat{\nabla} f p  dm \right| \le K \int_{\H} |\nabla f| p  dm.
\end{align*}

For fixed $R>0$  and
$I \subseteq I_0=\{(i,j)| 1\le i \le l, 1\le j \le k_i\}$, we set
$\Omega_I^R = \{(z,t)\in \H | |z_{i,j}| < R, \forall (i,j) \in I, |z_{i,j}| > R, \forall (i,j) \in I^c \}$ and
$H_I = H_I(z,t): (v_{1,1}, \cdots, v_{l,k_l}) \longrightarrow (w_{1,1}v_{1,1}, \cdots, w_{l,k_l}v_{l,k_l})$, where
$w_{i,j}=1$ when $(i,j) \in I$ and $w_{i,j} = \frac{z_{i,j}}{\bar{z}_{i,j}}$ when $(i,j) \in I^c$ is a linear transformation for fixed
$(z,t) \in \Omega_I^R$. We can choose $R$ large enough such that $B \subseteq \Omega_{I_0}^R$.

From the property of the heat kernel:
\begin{align*}
x_{i,j} \frac{\partial}{\partial y_{i,j}} p(z,t) = y_{i,j} \frac{\partial}{\partial x_{i,j}} p(z,t),
\end{align*}
we have
\begin{align*}
\bar{z}_{i,j} (\hat{\X}_{i,j} + \imath \hat{\Y}_{i,j})p = z_{i,j} (\X_{i,j} - \imath \Y_{i,j})p
\end{align*}

An integration by parts gives
\begin{align*}
\int_{\H} \X_{i,j} f p dm = - \int_{\H} f \X_{i,j} p dm, \qquad
\int_{\H} \Y_{i,j} f p dm = - \int_{\H} f \Y_{i,j} p dm, \\
\int_{\H} \hat{\X}_{i,j} f p dm = - \int_{\H} f \hat{\X}_{i,j} p dm, \qquad
\int_{\H} \hat{\Y}_{i,j} f p dm = - \int_{\H} f \hat{\Y}_{i,j} p dm.
\end{align*}

Note that $H_I$ has the property that
$\hat{\nabla}p = H_I \tilde{\nabla}_I p$ on $\Omega_I^R$,
where $\tilde{\nabla}_I =   (\tilde{\X}_{1, 1}, \tilde{\Y}_{1, 1}, \cdots, \tilde{\X}_{l, k_l}, \tilde{\Y}_{l, k_l})$,
$\tilde{\X}_{i,j} = \hat{\X}_{i,j}$, $\tilde{\Y}_{i,j} = \hat{\Y}_{i,j}$ for $(i,j) \in I$  and
$\tilde{\X}_{i,j} = \X_{i,j}$, $\tilde{\Y}_{i,j} = -\Y_{i,j}$ for $(i,j) \in I^c$.
One can refer to the proof in \cite{BBBC08} for more details.
If $f \in C_o^\infty(\H)$ and supported in $\Omega_I^R$ for some
$I$, we have
\begin{align*}
\int_{\H} \hat{\nabla}f p dm = -\int_{\H} f  \hat{\nabla}p dm = - \int_{\H} f H_I \tilde{\nabla}_I p dm.
\end{align*}
We divide into two parts. For $(i,j) \in I$, we have
\begin{align*}
\left| \int_{\H} f (\hat{\X}_{i,j} + \imath \hat{\Y}_{i,j})p dm \right|
&\le \left| \int_{\H} f (\hat{\X}_{i,j} + \imath \hat{\Y}_{i,j} - \X_{i,j} - \imath \Y_{i,j})p dm \right|
+\left| \int_{\H} f (\X_{i,j} + \imath \Y_{i,j})p dm \right| \\
& \le \left| \int_{\H} f (-4 a_i y_{i,j} + 4 \imath a_i x_{i,j})\partial_t p dm \right|
+\left| \int_{\H} (\X_{i,j} + \imath \Y_{i,j})f p dm \right| \\
& \le 4R \int_{\H} |f| |\partial_t p| dm + \int_{\H} |\nabla f | p dm  \\
& \lesssim_R \int_{\H} |\nabla f | p dm,
\end{align*}
where we use the assumption $m_f = 0$ (noticing that for $I \ne I_0$ we have $m_f = 0$ automatically), Lemma \ref{l2} and Theorem \ref{t2} in the last ``$\lesssim_R$''.
For the other case $(i,j) \in I^c$, we have
\begin{align*}
\left| \int_{\H} f (\hat{\X}_{i,j} + \imath \hat{\Y}_{i,j})p dm \right|
&\le \left| \int_{\H} f \frac{z_{i,j}}{\bar{z}_{i,j}} (\X_{i,j} - \imath \Y_{i,j})p dm \right| \\
& \le\left| \int_{\H} (\X_{i,j} - \imath \Y_{i,j})f \frac{z_{i,j}}{\bar{z}_{i,j}} p dm \right|
+\left| \int_{\H} f (\X_{i,j} - \imath \Y_{i,j})(\frac{z_{i,j}}{\bar{z}_{i,j}}) p dm \right| \\
& \lesssim  \int_{\H} |\nabla f | p dm
+ \left| \int_{\H} f \frac{1}{\bar{z}_{i,j}} p dm \right| \\
& \lesssim_R  \int_{\H} |\nabla f | p dm,
\end{align*}
where the last ``$\lesssim_R$'' follows from the assumption $m_f = 0$, Lemma \ref{l2} and Theorem \ref{t2}.
For the general case, we use partition of unity of the open cover $\mathop{\bigcup}\limits_{1 \le k \le n + 1, I \subset I_0} \Omega_I^{kR} $.
~ \hspace*{20pt} ~ \hfill $\Box$

\section*{Acknowledgement}
This work is partially supported by NSF of China (Grants No. 11625102). The author would like to thank M. Gordina for the useful bibliographic
information.

\nocite{*}
\bibliographystyle{abbrv}
\bibliography{gebib}

\mbox{}\\
Ye Zhang\\
School of Mathematical Sciences  \\
Fudan University \\
220 Handan Road  \\
Shanghai 200433  \\
People's Republic of China \\
E-Mail: 17110180012@fudan.edu.cn \quad or \quad zhangye0217@126.com \mbox{}\\

\end{document}